\documentclass{article}
\usepackage{epsfig}

\newtheorem{theorem}{Theorem}

\begin{document}

\title{On 3-manifolds} 

\author{Sergey Nikitin\\
Department of Mathematics\\
and\\
Statistics\\
Arizona State University\\
Tempe, AZ 85287-1804}

\maketitle 

 {\bf Abstract.} It is well known that a three dimensional (closed, connected and compact) manifold is obtained by identifying boundary faces from a polyhedron $P$. The study of $(\partial P)/\sim , $ the boundary $\partial P$ with the polygonal faces identified in pairs leads us to the following conclusion: either a three dimensional manifold is homeomorphic to a sphere or to a polyhedron $P$ with its boundary faces identified in pairs so that $(\partial P)/\sim$ is a two dimensional $CW$-complex. $(\partial P)/\sim$ is a finite number of two dimensional cells attached to each other along the edges of a finite graph that contains at least one closed circuit.  Each of those cells is obtained from a polygon where each side may be identified with one or more other different sides. Moreover, Euler characteristic of  $(\partial P)/\sim $ is equal to one and the fundamental group of $(\partial P)/\sim $ is not trivial.

\vspace{0.1cm}

\section{Main result}

\vspace{0.1cm}

It is well known \cite{SeifThrelf} that a closed, compact and connected three dimensional manifold $M$ is homeomorphic to a polyhedron $P$ with its boundary faces (polygonal faces which are parts of planes) identified in pairs. We follow \cite{SeifThrelf} and assume that each edge (straight-line segment) from the boundary $\partial P$ of $P$ is adjacent to exactly two boundary faces. The boundary $\partial P$ with the corresponding polygonal faces identified is denoted by $(\partial P)/\sim.$\\
Consider involution $p_0$ (permutation) defined on the faces of $P$ that maps each face into its peer. Then each edge from $(\partial P)/\sim$ corresponds to the class of edges from $\partial P$ that are glued together in $(\partial P)/\sim.$ Without loss of generality we can assume that each face $L$ of $P$ can not have two different edges that are glued together in $(\partial P)/\sim.$ Let us consider an edge $\alpha \in (\partial P)/\sim.$ Then one can introduce involution $p_\alpha $ on the set of faces from $ P.$ If a face $L$ does not have edges from $\alpha $ then
$$
p_\alpha (L) = L.
$$
Otherwise, there exists exactly one face $K$ that has the same edge from $\alpha$ as $L$ and we declare
$$
p_\alpha (L) = K.
$$
Now consider the group $G$ defined by generators $\{p_0,\;\{p_\alpha \}_\alpha \},$ where the set of involutions $\{p_\alpha \}_\alpha $ corresponds to the set of all edges from $(\partial P)/\sim .$ 

The group $G$ allows us to introduce the following important new concept.
The order $order(\alpha)$ is called the smallest integer $m>0$ such that
$$
(p_0 \cdot p_\alpha )^m = 1.
$$
The degree $deg((\partial P)/\sim)$ is defined as
$$
deg((\partial P)/\sim) = \max_\alpha \{order(\alpha) \}.
$$
In turn,
$$
deg(M) = \min_P \{ deg((\partial P)/\sim) \} ,
$$
where the minimum is taken over all polyhedra $P$ such that $M$ is homeomorphic to $P$ with the faces identified in pairs.
We call $M$ flat if $deg(M)=2.$ All flat three dimensional manifolds admit simple classification.
\begin{theorem}
\label{involution}
If $deg(M)=2$ then $M$ is homeomorphic to a polyhedron $P$ with boundary faces identified in pairs so that $(\partial P)/\sim$ is either homeomorphic to a disk or to projective plane. 
\end{theorem}
{\bf Proof.}
The boundary $\partial P$ of the polyhedron $P$ is homeomorphic to two dimensional sphere $S.$ The image of $(\partial P)/\sim$ under this homeomorphism is denoted by  $S/\simeq .$ 
If $M$ is a flat manifold  then the involution $p_0$ corresponds to a piecewise linear one-to-one involution on the two dimensional sphere $S.$  It follows from  \cite{Lopez} that there exists a circle $S^1$ that splits $S$ into two parts $S_+$ and $S_-$ such that
$$
p_0:\; S^1 \to S^1
$$
and $p_0$ is a one-to-one mapping that interchanges $S_+$ and $S_-.$
It is known \cite{SeifThrelf} that the Euler characteristic
$$
\chi (S/\simeq) =1.
$$
Hence, $S/\simeq $ is either a disk or a projective plane (see, e.g., \cite{Massey} for details).
Q.E.D.\\
 The classification of flat manifolds allows to tackle certain difficult problems. In particular, the next statement proves the famous Poincare conjecture for flat manifolds.
 
\begin{theorem} 
\label{Poincare}
If $M$ is compact, closed, connected and simply connected $3$-dimensional flat manifold then $M$ is homeomorphic to a sphere. 
\end{theorem}
{\bf Proof.}
If $M$ is a closed, connected, compact, flat $3$-dimensional manifold and $\pi(M)=\{1\}$ then by Theorem \ref{involution}  $(\partial P)/\simeq$ is homeomorphic to a disk.  Hence, $M$ is homeomorphic to three dimensional sphere (see \cite{Glaser}, \cite{Whitehead} for details).
Q.E.D.\\

Now we consider three dimensional manifolds that are not flat. Examples of such non-flat manifolds are provided by lens spaces.
  A lens shell $\ell (q,p)$ for coprime integers 
$$
q > p \ge 1
$$
can be defined as follows. Consider the function of a complex variable $z\in {\rm C},$
$$
p_0(z)=\left\{ \begin{array}{cc}
               \frac{e^{2\cdot  \pi \frac{p}{q} i }}{\bar z },\;\;&\mbox{ for } z\cdot \bar z \le 1,\\
               \frac{e^{-2\cdot  \pi \frac{p}{q} i}}{\bar z },\;\;&\mbox{ for } z\cdot \bar z  > 1,
               \end{array}
	\right.
$$
Now define the equivalence on $\{ {\rm C}, \infty \} $ as
$$
0 \sim \infty,\;\;\mbox{ and } z_1 \sim z_2 \;\;\mbox{as long as } z_1 = p_0^m(z_2)
$$
for some integer $m,$ where $ p_0^2(z)$ denotes $ p_0( p_0(z))$ and $p_0^m(z) = p_0(p_0^{m-1}(z)).$

$\{ {\rm C}, \infty \} $ is identified with the $2$-dimensional sphere $S$ (e.g. with stereographic projection). The topological space
$$
S/\sim
$$
is called lens shell $\ell (q,p).$ Notice that in the literature the lens space is defined as polyhedron $P$ with faces identified in pairs so that $(\partial P)/\sim$ is homeomorphic to $\ell (q,p).$ It is not difficult to show that 
$$
\pi(\ell (q,p) ) = {\rm Z_q}
$$
The points 
$$
e^{2\cdot m \pi \frac{p}{q} i }\;\;m=0,\;1,\;\dots q-1
$$
are all mapped into a single vertex of $S/\sim .$ There is only one edge $\alpha $ in $S/\sim $ and
$$
p_\alpha (z) = \frac{1}{\bar z}
$$
It is easy to see that
$$
(p_\alpha \cdot p_0)^q = 1
$$
In other words, lens shell $\ell(q,p)$ is not flat for $q>2.$\\
We call an edge $\alpha \in (\partial P)/\sim$ collapsible if there exists a face $F\in \partial P$ such that 
$$
p_\alpha (p_0 ( F )) = F.
$$
Consider the subgroup $G_2$ of group $G$ generated by
$$
\{p_\alpha :\;order(\alpha )= 2, \; \alpha \in (\partial P)/\sim , \;\alpha  \mbox{ is not collapsible } \}.
$$
Then the set of the polygonal faces of $\partial P$ is partitioned into the orbits of $G_2 .$  Those orbits are identified in pairs in order to obtain $(\partial P)/\sim .$ A pair of orbits identified in  $(\partial P)/\sim $ is obtained from a polygon with some of its sides (if not all) glued together, each side may either stay alone or to be attached to one or more different other sides. \\

We call edge $\alpha \in  (\partial P)/\sim $ non-flat if $order(\alpha )>2.$ Hence $(\partial P)/\sim $ is a two dimensional $CW$-complex (see \cite{Whitehead}) where each two dimensional cell ($2$-cell) corresponds to a pair of orbits of $G_2$ identified in  $(\partial P)/\sim .$ All $2$-cells are attached to a finite graph formed by non-flat edges from $(\partial P)/\sim .$

 Notice that the graph formed by non-flat edges has closed circuits when $2$-cell from $(\partial P)/\sim$ contains either a lens shell or a compact connected two dimensional manifold of non-trivial genus.\\
\begin{theorem} A three dimensional manifold is either homeomorphic to a sphere  or to a polyhedron $P$ with boundary faces identified in pairs so that $(\partial P)/\sim$ is a two dimensional $CW$-complex with finite number of $2$-cells attached to each other along edges of a finite graph that contains at least one closed circuit. Moreover, Euler characteristic $\chi ((\partial P)/\sim)=1.$
\end{theorem}
{\bf Proof.} The statement of this theorem was already proved for flat manifolds (see Theorem \ref{involution}).\\
 Now we consider non-flat manifolds.  If $M$ is not flat then it is homeomorphic to a polyhedron $P$ with faces identified in pairs and we can consider the graph $\Gamma $ that consists of all edges $\alpha \in (\partial P)/\sim $ such that $order(\alpha ) >2. $ 
If $deg(M)>2$ then  $\Gamma $ is not empty. Moreover, $\Gamma $ has at least one closed circuit. Indeed, if $\Gamma $ was a set of trees then one could consider it as a part of one single spanning tree. That means $(\partial P)/\sim$ consists of a finite number of $2$-cells attached to that tree. Each $2$-cell is  collapsible. Since $\Gamma $  is a set of trees then $(\partial P)/\sim$ neither contains shell spaces nor $2$-manifolds with non-trivial genus. On the other hand, it is known \cite{SeifThrelf} that $\chi ((\partial P)/\sim) =1.$ That means $(\partial P)/\sim$ is collapsible, and therefore,  $M$ is homeomorphic to three dimensional sphere (see \cite{Glaser}, \cite{Whitehead} for details).  That contradicts to the assumption that $deg(M )>2.$ Thus, $ deg(M)>2$ implies the existence of at least one closed circuit in $\Gamma .$\\
Q.E.D.\\
We complete this paper with the proof of the famous Poincar\'e conjecture.
\begin{theorem}
A three dimensional manifold is either homeomorphic to a sphere or to a polyhedron $P$ with its boundary faces identified in pairs so that $(\partial P)/\sim$ has non-trivial fundamental group.
\end{theorem}
{\bf Proof.} The statement of this theorem was already proved for flat manifolds (see Theorem \ref{Poincare}).\\ 
Consider a non flat manifold $M, \;\; deg(M)>2.$ $M$ is homeomorphic to a three-dimensional ball such that the boundary of this ball is a triangulated sphere $S$ with its $2$-simplexes identified in pairs. $S/\sim$ denotes $S$ with the corresponding $2$-simplexes identified. Without loss of generality we can assume that $S/\sim$ does not have collapsible $2$-simplexes. The graph $\Gamma $ formed by non flat edges of $S/\sim $ has at least one closed circuit. Take a minimal spanning tree for $\Gamma $ (the minimal tree in $S/\sim $ that contains all vertices from $\Gamma $ and the maximal possible number of arcs from $\Gamma $) and contract it to a point, say $x_0.$ This operation corresponds to the deformations of the triangulation on $S$ obtained as a sequence of the elementary steps $\{s_\beta \}_\beta .$ Each step $s_\beta $ corresponds to the edge $\beta $ of a $2$-simplex from $S$ and the deformation $s_\beta $ is depicted in Fig.1. 
\begin{figure}[tb]
  \begin{center}
       \psfig{file=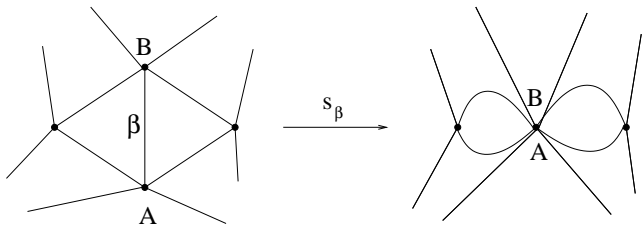,scale=0.4}         
         \caption{Deformation (and/or epimorphism) $s_\beta .$}
     \label{projectivePlane}
 \end{center}
\end{figure}
Consider a deformation $\bar s_\alpha $  for a non-flat edge $\alpha$ from $\Gamma .$  The result of this deformation $\bar s_\alpha (S/\sim )$ is obtained by identifying the corresponding $2$-cells from the following deformation for the triangulation of $S,$
$$
\prod_{\beta \in \alpha } s_\beta (S),
$$
where $\prod_{\beta \in \alpha } s_\beta$ denotes the superposition of deformations $s_\beta $ collected for all edges $\beta $ that are glued together in $S/\sim $ as the edge $ \alpha.$ Let $Tr(\Gamma )$ denote the spanning tree for $\Gamma $ that will be contracted to $x_0 .$ Then the deformation
$$
\prod_{\alpha \in Tr(\Gamma ) } \bar s_\alpha  (S/\sim )
$$
is obtained by identifying the corresponding $2$-cells from
$$
\prod_{\alpha \in Tr(\Gamma ) } \prod_{\beta \in \alpha } s_\beta (S).
$$
Due to its construction (Fig.1) each deformation $s_\beta $ can be identified with an epimorphism (we use the same notation $s_\beta $ for it) which is one-to-one outside the $1$-simplex with vertices $A,\;B$ depicted in Fig.1. Hence, we obtain the epimorphism of the fundamental group $\pi(S/\sim )$ onto the fundamental group
$$
\pi (\prod_{\alpha \in Tr(\Gamma ) } \bar s_\alpha  (S/\sim )).
$$
Consider $\prod_{\alpha \in Tr(\Gamma ) } \bar s_\alpha  (S/\sim ).$ The graph $\Gamma $ becomes a finite sum of circles $ \bigvee_{j=1}^{n+1} S^1_j $ having the only common point $x_0.$ Since $ \bigvee_{j=1}^{n+1} S^1_j $  comes from $\Gamma $ we call $\{ S^1_j \}_{j=1}^{n+1}$ non-flat circles.\\
The Euler characteristic for $\prod_{\alpha \in Tr(\Gamma ) } \bar s_\alpha  (S/\sim )$ is equal to one. Hence, if the number of non-flat circles is equal to one, then $\prod_{\alpha \in Tr(\Gamma ) } \bar s_\alpha  (S/\sim )$ is obtained from the sphere
$$
\prod_{\alpha \in Tr(\Gamma ) } \prod_{\beta \in \alpha } s_\beta (S)
$$
partitioned into exactly two flat pieces. Those two flat pieces has to be identified so that there is only one $2$-cell and one non-flat circle. It is possible to do only when $\prod_{\alpha \in Tr(\Gamma ) } \bar s_\alpha  (S/\sim )$  is a lens shell.\\
Indeed, consider
$$
S=\{z\cdot \bar z+x^2=1;\;\;(z,x)\in {\rm C}\times {\rm R} \},
$$
partitioned into two cells
$$
S_{+}=\{(z,x)\in S;\;x\ge 0\}, \;\; S_{-}=\{(z,x) \in S;\;x\le 0\}.
$$
Let $P$ denote the projection
$$
P(z,x)=z.
$$
Consider a homeomorphism
$$
g:\;\;S_{+}\;\to \; S_{-}
$$
such that after gluing $S_{+}$ with $S_{-}$ in accordance with $x \sim g(x)$ we obtain $S/\sim$ that has only one $2$-cell and only one non-flat circle. The mapping
$$
\nu (z)= P(g(z,\sqrt{1 -z\cdot \bar z}))
$$
is a homeomorphism of disk $D_2=\{z\cdot \bar z\le 1 \}$ to itself. Assume that the degree of the non-flat circle is equal to $q.$ Then, without loss of generality, we can consider only $\nu (z)$ such that there exists a natural number $p$ which is coprime with $q$ and
$$
\nu (z) = e^{i\frac{2\pi p}{q} } z \mbox{ for } \;\;z\cdot \bar z =1.
$$
It is well known \cite{Hempel}, \cite{Glaser}, that any self-homeomorphism of a disk which fixes the boundary is isotopic to the identity through self-homeomorphisms fixing the boundary. Hence, there exists an isotopy
$$
\mu;\;[0,1]\times D_2\;\to\;D_2
$$
such that
$$
\mu(0,z)=z\;\;\mbox{ and }\;\;e^{-i\frac{2\pi \cdot p }{q}}\nu (z)=\mu(1,z)\;\;\mbox{ for }\;z\cdot \bar z\le 1.
$$
Introduce notations
$$
\xi^t (z,x) = (e^{i\frac{2\pi \cdot p }{q}}\cdot \mu(t,z),\; -\sqrt{1-\mu(t,z)\cdot \overline{\mu(t,z)}}\;\;).
$$
$S^t/\sim$ denotes a complex which is obtained by gluing  $S_{+}$ and $S_{-}$ in accordance with $(z,x)\sim \xi^t (z,x).$
Due to its construction, $S^t/\sim$ is homeomorphic to $S^0/\sim$  for any $0\le t \le 1.$ On the other hand, $S^0/\sim$ is a lens shell and $S^1/\sim $ is $S/\sim .$

Hence, $\prod_{\alpha \in Tr(\Gamma ) } \bar s_\alpha  (S/\sim )$ is a lens shell. Thus, there exists an epimorphism of the fundamental group for  $S/\sim $ onto the fundamental group of the lens shell. That yields non-triviality of the fundamental group for $S/\sim $ under the condition that there is only one non-flat circle in $\prod_{\alpha \in Tr(\Gamma ) } \bar s_\alpha  (S/\sim ).$ \\
Introduce notations
$$
epi(S/\sim )= \prod_{\alpha \in Tr(\Gamma ) } \bar s_\alpha  (S/\sim )
$$
and $epi(S)$ denotes the corresponding sphere. That means $epi(S/\sim )$ is obtained from $epi(S)$ by identifying in pairs the corresponding $2$-cells. \\
We conduct the proof by redactio ad absurdum. Assume 
$$
\pi (S/\sim)=\{1\}.
$$
Then
$$
\pi (epi(S/\sim ))=\{1\}.
$$
Hence, $S^1_{n+1}$ is homotopic to $x_0.$ That yields the existence of a continuous mapping 
$$
\varphi :\;\;\{(r\cos(\theta ),\;r\sin(\theta ));\;0\le r \le 1,\;\;0\le \theta < 2\pi \} \;\;\to \;\; epi(S/\sim )
$$
such that  
$$
\varphi (0, \theta ) = x_0\;\;\;\mbox{for}\;\;\; \theta \in [0,2\pi]
$$
and $\{\varphi (1, \theta );\;0\le \theta< 2\pi \}$ is a parameterization of $S_{n+1} \subset epi(S/\sim).$ Without loss of generality one can assume that
$$
\varphi (r,0)=\varphi (r,2\pi)=x_0\;\;\;\mbox{ for }\;\;0\le r \le 1
$$
and
$$
\varphi (1, \theta_1)\not= \varphi (1, \theta_2)\;\;\mbox{ for }\;\;0\le \theta_1< \theta_2<2\pi.
$$
There are only two possibilities:
\begin{itemize}
\item[(i)] $epi(S/\sim)\setminus \varphi(D_2) \not= \emptyset , $
\item[(ii)] $epi(S/\sim) = \varphi(D_2).$
\end{itemize}
Consider (i). Then there exists a continuous epimorphism
$$
\beta_1:\;\;epi(S/\sim)\;\;\to\;\;\beta_1(epi(S/\sim)),
$$
such that $\beta_1$ maps $\varphi(D_2)$ to $x_0$ and  $\beta_1$  is one-to-one on 
$$
epi(S/\sim)\setminus \varphi(D_2) .
$$
Notice that Euler characteristic
$$
\chi (\beta_1(epi(S/\sim)))=1.
$$
Hence, $\beta_1(epi(S/\sim))$ has at least one non-flat circle. Moreover, $\beta_1(epi(S/\sim))$ is obtained from a sphere by identifying in pairs flat $2$-cells.\\

Consider (ii). The two-dimensional complex $epi(S/\sim)$ can be written as
$$
epi(S/\sim) = \cup_{j=1}^{n+1}F_j,
$$
where each $F_j$  corresponds to a pair of $2$-cells from $epi(S).$ Then
$$
D_2=\cup_{j=1}^{n+1} \varphi^{-1}(F_j)
$$
and
$$
int(\varphi^{-1}(F_j)) \cap int(\varphi^{-1}(F_i)) = \emptyset\;\;\mbox{ for } \;\;i\not=j,
$$
where $int(Q)$ denotes the interior of $Q.$ Moreover, for their closures we have
$$
\overline{\varphi^{-1}(F_j)} \cap \overline{\varphi^{-1}(F_i)} \subset  \varphi^{-1}(\bigvee_{j=1}^{n+1} S^1_j) ,
$$
where $i\not= j .$ On the other hand,
$$
\varphi^{-1}(S^1_{i}) \cap \varphi^{-1}(S^1_j))\subset \varphi^{-1}(x_0)
$$
where $i\not= j .$ That means the disk $D_2$ is partitioned in the way similar to the one depicted in Fig.2.
\begin{figure}[tb]
  \begin{center}
           \psfig{file=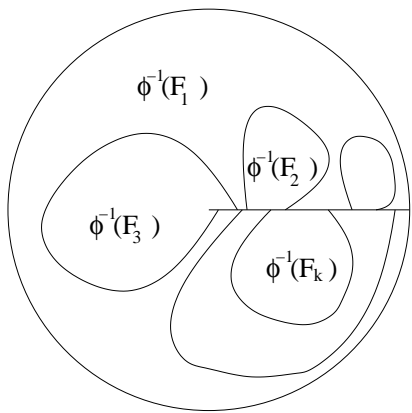,scale=0.4}         
            \caption{Partition of the disk.}
            \label{projectivePlane}
         \end{center}
\end{figure}
The interior of each connected component of $\varphi^{-1}(F_j)$ is simply-connected $(j=1,2, \dots n+1 ).$  $D_2$ is compact. That yields the existence $\overline{\varphi^{-1}(F_k)},$ such that $\partial F_k=S^1_j,$ where  $S^1_j$ is a non-flat circle. Hence, there exists an epimorphism 
$$
\beta_1: epi(S/\sim )\;\to \;\beta_1(epi(S/\sim ))
$$
such that $\beta_1$ maps $\overline{F_k}$ to $x_0$ and $\beta_1$ is one-to-one on
$$
epi(S/\sim ) \setminus \overline{F_k} .
$$
By construction, $\beta_1(epi(S/\sim ))$ has non-flat circles and $\beta_1(epi(S/\sim ))$ is obtained by identifying in pairs $2$-cells from a partition of sphere.\\

After a finite number of steps we obtain a continuous epimorphism
$$
\gamma = \beta_k\circ \beta_{k-1}\circ \dots \beta_1 ,
$$
where $\beta_k\circ \beta_{k-1}\circ \dots \beta_1$ denotes superposition of continuous epimorphisms 
$$
\beta_k,\; \beta_{k-1}, \dots \beta_1.
$$
Each $\beta_j\;\;(j=1,2,\dots k)$ is constructed in the same way as $\beta_1.$ Only $\beta_j$ corresponds to a different non-flat circle which meets one of the conditions (i), (ii).
$$
\gamma:\;\;epi(S/\sim )\; \to \; \gamma (epi(S/\sim ))
$$
and $\gamma (epi(S/\sim ))$ has the only non-flat circle. Moreover, $\gamma (epi(S/\sim ))$ is obtained by identifying $2$-cells from a sphere partitioned into exactly two flat pieces. Hence, 
$$
\gamma (epi(S/\sim ))
$$
is a lens shell and
$$
\pi (\gamma (epi(S/\sim ))) \not= \{1\}.
$$
However, it contradicts to the assumption that $\pi(S/\sim)$ is trivial. This contradiction completes the proof.
Q.E.D.\\

\bibliographystyle{ams-plain}

\end{document}